\documentclass[12pt,a4paper,reqno]{amsart}
\usepackage{amsmath}
\usepackage{amsthm}
\usepackage{amsfonts}
\usepackage{amssymb}
\usepackage{setspace}
\numberwithin{equation}{section}

\theoremstyle{plain}

\newtheorem{theorem}[subsection]{Theorem}
\newtheorem{proposition}[subsection]{Proposition}
\newtheorem{lemma}[subsection]{Lemma}
\newtheorem{corollary}[subsection]{Corollary}

\theoremstyle{definition}

\renewcommand{\leq}{\leqslant}
\renewcommand{\geq}{\geqslant}

\newsavebox{\proofbox}
\savebox{\proofbox}{\begin{picture}(7,7)%
 \put(0,0){\framebox(7,7){}}\end{picture}}

\def\F{\mathbb{F}}

\DeclareMathOperator{\Alg}{Alg}
\DeclareMathOperator{\Spec}{Spec}

\def\proof{\noindent\textit{Proof. }}
\def\endproof{\hfill{\usebox{\proofbox}}\vspace{11pt}}

\begin{document}

\onehalfspace
\title{Sum-product phenomena in $\F_p$: a brief introduction}

\author{Ben Green}
\address{Centre for Mathematical Sciences\\
Wilberforce Road\\
     Cambridge CB3 0WA\\
     England
}
\email{b.j.green@dpmms.cam.ac.uk}

\thanks{The author holds a Leverhulme Prize and is grateful to the Leverhulme Trust for their support.}

\begin{abstract}
These notes arose from my Cambridge Part III course on Additive Combinatorics, given in Lent Term 2009. The aim was to understand the simplest proof of the Bourgain-Glibichuk-Konyagin bounds for exponential sums over subgroups. As a byproduct one obtains a clean proof of the Bourgain-Katz-Tao theorem on the sum-product phenomenon in $\F_p$. The arguments are essentially extracted from Bourgain's paper \cite{bourgain-multi}, and I benefitted very much from being in receipt of unpublished course notes of Elon Lindenstrauss. No originality is claimed. 
\end{abstract}
\maketitle

\section{introduction}

If $A$ and $B$ are subsets of a ring then we write $A + B := \{a + b : a \in A, b \in B\}$ and $A \cdot B := \{ab : a \in A, b \in B\}$.

The Bourgain-Katz-Tao theorem \cite{bkt-paper} is a celebrated result roughly stating that there are no approximate subrings of $\F_p$. It has found a great many applications.

\begin{theorem}[Bourgain-Katz-Tao]\label{bkt}
Suppose that $\delta > 0$ and that $A \subseteq \F_p$ is a set of cardinality between $p^{\delta}$ and $p^{1 - \delta}$. Then there is some $\delta' > 0$ such that either $A \cdot A$ or $A + A$ has cardinality $\gg |A|^{1 + \delta'}$.
\end{theorem}

The requirement that $|A| \geq p^{\delta}$ was later dropped in work of Bourgain, Glibichuk and Konyagin.  Shortly after the appearance of the Bourgain-Katz-Tao result, the following theorem was proved by Bourgain, Glibichuk and Konyagin.

\begin{theorem}[Bourgain-Gilibichuk-Konyagin]\label{bgk-theorem}
Suppose that $H \leq \F_p^{\times}$ is a multiplicative subgroup of size at least $p^{\delta}$. Then uniformly in $\xi \neq 0$ we have
\[ \frac{1}{|H|} |\sum_{x \in H} e(x\xi/p)| \ll p^{-\delta'},\] where $\delta' = \delta'(\delta) > 0$.
\end{theorem}

This states that multiplicative subgroups of $\F_p^{\times}$, even very small ones, have very little additive structure indeed: they are ``additively pseudorandom''. Theorem \ref{bgk-theorem} is something of a triumph for additive combinatorics, for the question had previously been extensively studied by quite sophisticated number-theoretical arguments. The best result obtained with such methods applies only when $\delta > 1/4$; it is due to Konyagin \cite{konyagin}.

The first published proof of Theorem \ref{bgk-theorem} utilised the Bourgain-Katz-Tao theorem, but subsequently it has been realised that Theorems \ref{bkt} and the Bourgain-Katz-Tao theorem may be derived separately from some rather simpler considerations. We will partially follow a recent paper of Bourgain \cite{bourgain-multi} as well as some unpublished notes of Elon Lindenstrauss, and give what seems to be about the simplest derivation of these two theorems. As in the original version, we shall only establish the Bourgain-Katz-Tao theorem under the assumption that $|A| \geq p^{\delta}$, and furthermore our method for proving Theorem \ref{bgk-theorem} does not give the best known dependence of $\delta'$ on $\delta$ as $\delta \rightarrow 0$.

Both theorems are manifestations of the \emph{sum-product phenomenon}: additive and multiplicative structure find it hard to coexist.

\section{Rough notation, Ruzsa calculus and Balog-Szemer\'edi-Gowers}

One of my aims when giving the course \cite{green-course} was to introduce various types of rough notation which, once comprehended, seem to make arguments easier to follow.

The most basic convention in force is that the letters $c$ and $C$ stand for absolute constants with $0 < c < 1 < C$, but that different instances of the notation, even on the same line, might refer to different constants. For all of the arguments in these notes each instance of one of these letters could be turned into a concrete constant by a keen reader. Occasionally we will use subscripts to denote dependence on other parameters; for example $c_{\alpha,\beta}$ denotes a small positive absolute constant depending on $\alpha$ and $\beta$. \vspace{11pt}

\textsc{Rough notation at scale $K$.} Let $K \geq 2$ be a parameter. If $X$ and $Y$ are quantities then we write $X \lessapprox Y$ or $Y \gtrapprox X$ to mean $X \leq K^CY$, and $X \approx Y$ if $X \lessapprox Y$ and $Y \lessapprox X$. Different instances of the notation are allowed to involve different implied absolute constants $C$. The parameter $K$ measures the ``roughness'' of this notation and it will be fixed for the duration of any given argument and referred to as the \emph{roughness parameter}.\vspace{11pt}

\textsc{Ruzsa calculus.} This is a notation for handling a number of inequalities for the size of sumsets discovered by Imre Ruzsa. Let $K \geq 2$ be a roughness parameter. If $A$ and $B$ are finite sets in some abelian group then we write $A \sim B$ if $|A - B| \lessapprox |A|^{1/2}|B|^{1/2}$ in rough notation at scale $K$. Equivalently, $A \sim B$ if the \emph{Ruzsa distance} between $A$ and $B$, as introduced in \cite[p. 60]{tao-vu}, is $O(\log K)$. In the course \cite{green-course} the following collection of statements was referred to as ``Ruzsa calculus''. In this proposition we write $\sigma[A] := |A+A|/|A|$ for the \emph{doubling constant} of the set $A$. 

\begin{proposition}[Ruzsa Calculus] Suppose that $U,V$ and $W$ are sets in some ambient abelian group. We use rough notation at some fixed scale $K$.\begin{enumerate}
\item Suppose that $U \sim V$. Then $U \sim -V$, $|U| \approx |V|$ and $\sigma[U], \sigma[V] \approx 1$. 
\item If $U \sim V$ and $V \sim W$, then $U \sim W$.
\item Suppose that $U \sim V$ , that $\sigma[W] \approx 1$ and that there is some $x$ such that $|U \cap (x + W)| \approx |U| \approx |W|$. Then $U \sim V \sim W$.  
\item Suppose that $\sigma[U],\sigma[W] \approx 1$ and that there is some $x$ such that $|U \cap (x + W)| \approx |U| \approx |W|$. Then $U \sim W$.
\item Suppose that $U \sim V \sim W$. Then $U \sim V + W$.\end{enumerate}
\end{proposition}
Proofs of these statements may be found in \cite[Chapter 2]{tao-vu} or in the second set of notes from my course \cite{green-course}.\vspace{11pt}

\textsc{Additive energy and Balog-Szemer\'edi-Gowers}. If $A$ and $B$ are two finite sets in some abelian group then we write $\omega_+(A,B)$ for $|A|^{-3/2}|B|^{-3/2}$ times the number of quadruples $(a_1,b_1,a_2,b_2) \in A \times B \times A \times B$ with $a_1 + b_1 = a_2 + b_2$. This quantity, called the additive energy, is easily seen to lie in the interval $[0,1]$. The property of having large additive energy is an important notion of structure because (as we shall see) it often arises from Fourier-analytic arguments. Remarkably it is closely related to the more precise notion of small doubling, and one then has the tools of Ruzsa calculus at one's disposal.

\begin{proposition}[Balog-Szemer\'edi-Gowers]\label{bsg}
Let $K \geq 2$ be a roughness parameter and suppose that $A,B$ are finite subsets of some abelian group. Then
\begin{enumerate}
\item If $A \sim B$ then $\omega_+(A,B) \approx 1$.
\item If $\omega_+(A,B) \approx 1$ then there are subsets $A' \subseteq A$ and $B' \subseteq B$ with $|A'| \approx |A|$ and $|B'| \approx |B|$ such that $A' \sim B'$.
\end{enumerate}
\end{proposition}

Statement (i) is rather easy and is nothing more than a single application of the Cauchy-Schwarz inequality: the Balog-Szemer\'edi-Gowers theorem is (ii). It is absolutely necessary to pass to subsets $A'$ and $B'$ in general. This is because the property of having large additive energy persists under additing a few arbitrary elements to the sets involved, whereas the property of having small sumset does not. For example if $A = \{1,\dots,n\} \cup \{2^n, 2^{2n},\dots, 2^{n^2}\}$ then $\omega_+(A,A) \approx 1$ but $|A-A|$ has size $c|A|^2$.

The proof of the Balog-Szemer\'edi-Gowers theorem in this form may be found in \cite{tao-vu} or in the fourth set of notes from my course \cite{green-course}.

\section{Simple instances of the sum-product phenomenon}

The Bourgain-Katz-Tao theorem asserts that a set $A$ must grow under either addition or multiplication. It is rather easier to establish that $A$ grows under a small number of additions \emph{and} multiplications, and for many applications results of this type suffice. Perhaps the simplest such result is due to Glibichuk and Konyagin. Before establishing it we prove a lemma.

\begin{lemma}\label{lemma6.1}
Suppose that $A \subseteq \F_p$ is a set. Then there is some $\xi \in \F^{\times}_p$ such that $|A + \xi \cdot A| \geq \frac{1}{2}\min(|A|^2,p)$. 
\end{lemma}
\noindent\emph{Remark.} Here we have written $\xi \cdot A := \{\xi a : a \in A\}$.\vspace{11pt}

\proof There is one key idea, which is to use additive energy. The sum $\sum_{\xi \neq 0} \omega_+(A, \xi \cdot A)$ counts $|A|^{-3}$ times the number of solutions to $a_1 - a_2 = \xi(a_3 - a_4)$. For each of the $|A|^2(|A|-1)^2$ quadruples $(a_1,a_2,a_3,a_4)$ with $a_1 \neq a_2$ and $a_3 \neq a_4$ there is a unique choice of $\xi \neq 0$ which satisfies this equation. If $a_1 = a_2$ and $a_3 = a_4$ then any choice of $\xi$ works, and hence
\[ \sum_{\xi \neq 0} \omega_+(A,\xi \cdot A) = \frac{1}{|A|}\big(p - 1+ (|A|-1)^2\big).\]
In particular there is some $\xi$ for which 
\[ \omega_+(A,\xi \cdot A) \leq \frac{1}{|A|} + \frac{(|A| - 1)^2}{|A|p} \leq \frac{1}{|A|} + \frac{|A|}{p} \leq 2\max(\frac{1}{|A|}, \frac{|A|}{p}).\] The result follows immediately from the (simple direction of) the relationship between additive energy and sumsets, that is to say Proposition \ref{bsg} (i).
\endproof

Here is the promised result of Glibichuk and Konyagin, asserting that we have growth under addition \emph{and} multiplication. Recall that if $k \geq 1$ is an integer we write $kA = A + \dots + A$ and $A^k = A \cdot A \cdot \dots \cdot A$, where in each case there are $k$ copies of the set $A$.

\begin{theorem}[Growth under addition and multiplication]\label{add-and-multiply}
Suppose that $A \subseteq \F_p$ is a set. Then we have $|3A^2 - 3A^2| \geq \frac{1}{2}\min(|A|^2,p)$.
\end{theorem}
\proof We begin with the observation that $|A + \xi \cdot A| = |A|^2$ unless there are distinct pairs $(a_1,a_2)$ and $(a_3,a_4)$ with $a_1 + \xi a_2 = a_3 + \xi a_4$, which means that $\xi \in \frac{A-A}{A-A}$. Suppose first of all that $\frac{A - A}{A - A}$ is not the whole of $\F_p$. Then there is some $\xi \in \frac{A-A}{A-A}$ for which $\xi + 1 \notin \frac{A-A}{A-A}$ which, by the preceding remarks, implies that $|A + (\xi + 1)\cdot A| = |A|^2$. Supposing that $\xi = \frac{a_1 - a_3}{a_2 - a_4}$, we have
\[ 3A^2 - 3A^2 \supseteq (a_2 - a_4) \cdot A + (a_1 - a_3 + a_2 - a_4) \cdot A \supseteq (a_2 - a_4) \cdot (A + (\xi +1)\cdot A).\] This implies that $|3A^2 - 3A^2| = |A|^2$.
Alternatively, suppose that $\frac{A-A}{A-A} = \F_p$. Then by Lemma \ref{lemma6.1} there is some $\xi$ such that $|A + \xi \cdot A| \geq \frac{1}{2}\min(|A|^2,p)$. Supposing that $\xi = \frac{a_1 - a_3}{a_2 - a_4}$ we may proceed much as before, observing now that
\[ 3A^2 - 3A^2 \supseteq 2A^2 - 2A^2 \supseteq (a_2 - a_4)\cdot (A + \xi \cdot A).\]
This concludes the proof.\endproof

By repeated application of this theorem it is not hard to obtain the following corollary.

\begin{corollary}[Sums and products generate the whole of $\F_p$]\label{cor6.1}
Let $\delta$, $0 < \delta < 1$, be a parameter and suppose that $A \subseteq \F_p$ is a set of cardinality at least $p^{\delta}$. Then there is some $k = k(\delta)$ such that $k A^k - kA^k = \F_p$. In particular if $A$ is a multiplicative subgroup of $\F^{\times}_p$ then $k A - kA = \F_p$.
\end{corollary}
\noindent\emph{Remark.} The constant $k(\delta)$ could be taken to be $e^{C/\delta^C}$ if desired.\vspace{11pt}

\proof By repeated application of Theorem \ref{add-and-multiply} we can certainly find a $k = k(\delta)$ such that $|kA^k - kA^k| \geq p/2$. The inequality is strict since $p > 2$. It remains only to note that if $X \subseteq \F_p$ has size greater than $p/2$ then $X + X = \F_p$: indeed for any $t \in \F_p$ the sets $X$ and $t - X$ must intersect or else their union would have size greater that $p$.\endproof

In the above discussion we considered sumsets $A + \xi \cdot A$. For the next few paragraphs only let us write $\Alg_K(A) = \{\xi \in \F^{\times}_p : |A + \xi \cdot A| \leq K|A|\}$. Noting that $\xi \in \Alg_{K^C}(A)$ if and only if $A \sim \xi \cdot A$, the following collection of properties follow easily by Ruzsa calculus.

\begin{lemma}[Properties of $\Alg$]\label{alg-properties} Let $K \geq 2$. We have the following statements.
\begin{enumerate}
\item If $\xi \in \Alg_K(A)$ then $\xi^{-1} \in \Alg_K(A)$.
\item If $\xi \in \Alg_K(A)$ then $-\xi \in \Alg_{K^C}(A)$.
\item If $\xi_1,\xi_2 \in \Alg_{K}(A)$ and $\xi_1 \neq -\xi_2$ then $\xi_1 + \xi_2 \in \Alg_{K^C}(A)$.
\item If $\xi_1,\xi_2 \in \Alg_K(A)$ then $\xi_1\xi_2 \in \Alg_{K^C}(A)$.
\end{enumerate}
\end{lemma}
\proof Only the last is not an absolutely immediate application of Ruzsa calculus. For it we need the additional observation that if $A \sim \xi_2\cdot A$ then $\xi_1 \cdot A \sim \xi_1 \xi_2 \cdot A$ for any $\xi_1 \neq 0$.\endproof

Using this, we can obtain the following rather pleasing expression of the sum-product phenomenon in $\F_p$.

\begin{corollary}\label{corollary6.2a}
Let $\alpha$ and $\beta$, $0 < \alpha,\beta < 1$, be parameters. Suppose that $A,B \subseteq \F_p$ are two sets with $p^{\alpha} \leq |A| \leq p^{1-\alpha}$ and $|B| \geq p^{\beta}$. Then there is $b \in B$ such that $|A + b\cdot A| \geq |A|^{1 + c_{\alpha,\beta}}$.
\end{corollary}
\proof Let $k = k(\beta)$ be the number whose existence is guaranteed by Corollary \ref{cor6.1}, so that $k B^k - kB^k = \F_p$. Suppose that $|A + b \cdot A| \leq K|A|$ for all $b \in B$, that is to say suppose $B \subseteq \Alg_K(A)$. Then by repeated applications of Lemma \ref{alg-properties} we have 
\[ \F^{\times}_p = (kB^k - kB^k) \setminus \{0\} \subseteq \Alg_{K^C}(A),\] where $C$ depends only on $k$ and hence, ultimately, on $\beta$. Thus
\[ |A + \xi \cdot A| \leq K^{C_{\beta}}|A|\] for all $\xi \in \F^{\times}_p$. This, however, may be contrasted with Lemma \ref{lemma6.1}, which guarantees a $\xi$ for which
\[ |A + \xi \cdot A| \geq \frac{1}{2}\min(|A|^2,p).\] The corollary is immediate.\endproof

\section{Additive-multiplicative Balog-Szemer\'edi-Gowers Theorem}

Corollary \ref{corollary6.2a} is particularly powerful in conjunction with the next result.

\begin{proposition}[Bourgain]\label{ambsg}
Suppose that $K \geq 2$ is a roughness parameter and that $A \subseteq \F_p$ and $B \subseteq \F^{\times}_p$ are sets. Suppose that $\omega_+(A,b\cdot A) \approx 1$ for all $b \in B$. Then there is an $x$ and subsets $A' \subseteq A$ and $B' \subseteq xB$ with $|A'| \approx |A|$ and $|B'| \approx |B|$ such that $A' \sim b' \cdot A'$ for all $b' \in B'$.
\end{proposition}
\proof From the Balog-Szemer\'edi-Gowers theorem for each $b \in B$ there are sets $X_b, Y_b \subseteq A$ with $|X_b|,|Y_b| \approx |A|$ and $|X_b + b \cdot Y_b| \approx |A|$. The task, of course, is to somehow remove the dependence of these sets on $b$.

We pause to observe a consequence of the Cauchy-Schwarz inequality which is useful in many situations.

\begin{lemma}\label{cbs-sets}
Suppose that $S$ is a finite set and that $S_1,\dots,S_k$ are subsets of $S$ with $|S_1|,\dots |S_k| \geq \delta$. Then there is some $i$ such that $|S_i \cap S_j| \geq \delta^2/2$ for at least $\delta^2 k/2$ values of $j$. 
\end{lemma}
\proof We have $\sum_{i=1}^k\sum_{x \in S} 1_{S_i}(x) = \delta k |S|$. By the Cauchy-Schwarz inequality this implies that $\sum_{x \in S} |\sum_{i =1}^k 1_{S_i}(x)|^2 \geq \delta^2 k^2 |S|$,
that is to say $\sum_{i,j} |S_i \cap S_j| \geq \delta^2 k^2 |S|$.
In particular there is some $i$ for which $\sum_j |S_i \cap S_j| \geq \delta^2 k |S|$, and this easily implies the lemma.\endproof

 Applying this with $S = A \times A$ and the sets $S_i$ equal to $X_b \times Y_b$, $b \in B$ we see that there is some $b_0 \in B$ and a set $B'$ with $|B'| \approx |B|$ such that 
\[ |X_{b_0} \cap X_b| , |Y_{b_0} \cap Y_b| \approx |A| \qquad \mbox{for all $b \in B'$}.\] Suppose that $b \in B'$. By assumption we have $X_b \sim b \cdot Y_b$ and $X_{b_0} \sim b_0 \cdot Y_{b_0}$. In particular $\sigma[X_{b_0}],\sigma[X_b],\sigma[Y_{b_0}],\sigma[Y_b] \approx 1$ and so by rule (iv) of Ruzsa calculus we have $X_{b_0} \sim X_b$ and $Y_{b_0} \sim Y_b$, the second of which implies $b \cdot Y_{b_0} \sim b \cdot Y_b$. Comparing all of these equivalences we see that 
\[ b_0 \cdot Y_{b_0} \sim X_{b_0} \sim X_b \sim b \cdot Y_b \sim b \cdot Y_{b_0},\] and hence $Y_{b_0} \sim (b/b_0)\cdot Y_{b_0}$. Redefining $B'$ to be the set of $b/b_0$ and setting $A' := Y_{b_0}$ concludes the proof of Proposition \ref{ambsg}.\endproof

Let us record the result of combining Proposition \ref{ambsg} with Corollary \ref{corollary6.2a}.

\begin{corollary}\label{corollary6.2}
Let $\alpha$ and $\beta$, $0 < \alpha,\beta < 1$, be parameters. Suppose that $A,B \subseteq \F_p$ are two sets with $p^{\alpha} \leq |A| \leq p^{1-\alpha}$ and $|B| \geq p^{\beta}$. Then there is $b \in B$ such that $\omega_+(A,b\cdot A) \geq |A|^{1 + c_{\alpha,\beta}}$.
\end{corollary}
\proof Suppose not: then there is some roughness parameter $K = |A|^{o(1)}$ such that $\omega_+(A,b\cdot A)\approx 1$ for all $b \in B$. Using Proposition \ref{ambsg} we may extract subsets $A'$ and $B'$ with $|A'| \geq p^{\alpha/2}$ and $|B'| \geq p^{\beta/2}$ such that $A' \sim b' \cdot A'$ for all $b' \in B'$. This, however, is contrary to Corollary \ref{corollary6.2}.\endproof

\section{The Bourgain-Katz-Tao theorem}

In this section we prove Theorem \ref{bkt}, namely the lower bound on $\min(|A + A|,|A \cdot A|)$ under the assumption that $p^{\delta} \leq |A| \leq p^{1 - \delta}$.  The key ingredient is the additive-multiplicative Balog-Szemer\'edi-Gowers theorem of Bourgain.

Suppose then that $|A + A|, |A\cdot A| \approx |A|$, where the underlying roughness parameter $K$ has size $p^{o(1)}$. 

Consider first the sets $a \cdot A$, $a \in A$. An application of Lemma \ref{cbs-sets} implies that there is some $a_0 \in A$ such that \[ |A \cap \frac{a_0}{a}\cdot A| = |a \cdot A \cap a_0\cdot A| \gtrapprox |A|\] for $\gtrapprox |A|$ values of $a \in A$.  Take such an $a$, and consider $A' := A \cap \frac{a}{a_0}A$.  Then the doubling $\sigma[A']$ is $\lessapprox 1$, and hence the additive energy $\omega_+(A',A')$ is $\approx 1$. It follows that $\omega_+(A, \frac{a}{a_0}A) \gtrapprox 1$, this bound being uniform in $a$. This is contrary to the additive-multiplicative Balog-Szemer\'edi-Gowers theorem, applied with $B$ equal to the set of these ratios $a/a_0$.\endproof

As we remarked earlier the requirement that $|A| \geq p^{\delta}$ is not necessary. It may be removed by combining Theorem \ref{add-and-multiply} with the \emph{Katz-Tao lemma}, one instance of which is the following statement: if $|A+A| \approx |A\cdot A| \approx |A|$ (using rough notation at some scale $K$) then there is a set $A' \subseteq A$ with $|A'| \approx |A|$ and for which $|3A^{\prime 2} - 3A^{\prime 2}| \approx |A'|$. We omit the proof, which involves a good deal of Ruzsa calculus as well as an application of the Balog-Szemer\'edi-Gowers theorem.

\section{Additive structure of the large spectrum}

We turn now to the proof of Theorem \ref{bgk-theorem}, the Bourgain-Glibichuk-Konyagin estimate for exponential sums or ``Gauss sums'' over multiplicative subgroups. 

Suppose that $A \subseteq \F_p$ is a set and that $\alpha$, $0 < \alpha < 1$,  is a parameter. Then we define
\[ \Spec_{\alpha}(A) := \{ \xi \in \F_p : |\sum_{x \in A} e(\xi x/p)| \geq \alpha |A|\}.\]
This is manifestly a symmetric set containing zero. 
If $A = H$ is a multiplicative subgroup of $\F_p$ then it is very easy to see that $\Spec_{\alpha}(A)$ is multiplicatively $H$-invariant: this follows from the fact that $1_H(hx) = 1_H(x)$ for all $h \in H$ and all $X \in \F_p$. 

If the theorem is not true then $\Spec_{\eta}(H)$ contains a nonzero element for some $\eta = p^{-\delta'}$.  The strategy is then to consider various $\Spec_{\alpha}(H)$, each of which will certainly be nonempty and $H$-invariant. The key additional observation is that $\Spec_{\alpha}(H)$ has a certain amount of weak \emph{additive} structure, too. We establish this fact, which does not depend on any multiplicative properties of $H$, in the following lemma. 

\begin{lemma}[Additive structure of Spec]\label{add-spec}
Suppose that $H \subseteq \F_p$ is any set and that $\alpha$, $0 < \alpha < 1$, is a parameter. Suppose that $B$ is a subset of $\Spec_{\alpha}(H)$. Then for a proportion at least $\frac{1}{2}\alpha^2$ of the pairs $x,y \in B$ we have $x - y \in \Spec_{\alpha^2/2}(H)$. 
\end{lemma}
\proof Write $B = \{\xi_1,\dots,\xi_k\}$. Then for each $i = 1,\dots,k$ there is a unit modulus complex number $c_i$ such that 
\[ c_i \sum_{x \in H} e(\xi_i x/p) \geq \alpha |H|.\]
Summing, we have
\[ \sum_{x \in H}\sum_{i=1}^k c_i e(\xi_i x/p) \geq k\alpha |H|.\] 
By the Cauchy-Schwarz inequality this implies that
\[ \sum_{x \in H}\sum_{i=1}^k\sum_{j=1}^k c_i \overline{c_j} e((\xi_i - \xi_j) x/p) \geq k^2 \alpha^2 |H|.\]
By the triangle inequality we thus have
\[ \sum_{i=1}^k \sum_{j=1}^k |\sum_{x \in H} e((\xi_i - \xi_j)x/p)| \geq k^2 \alpha^2 |H|.\]
Suppose that the conclusion of the lemma were false. Then, using the trivial bound $|\sum_{x \in H} e(\eta x/p)| \leq |H|$, the contribution to the left hand side from those pairs $(i,j)$ for which $\xi_i - \xi_j \in \Spec_{\alpha^2/2}(H)$ would be less than $\alpha^2 k^2|H|/2$. Furthermore the contribution from each pair $(i,j)$ for which $\xi_i - \xi_j \notin \Spec_{\alpha^2/2}$ is, by definition, no more than $\alpha^2|H|/2$, and so these terms contribute at most $\alpha^2 k^2 |H|/2$ in total. This contradiction establishes the lemma.\endproof

We will only use this lemma in the case when $B$ is the whole of $\Spec_{\alpha}(H)$, but for other applications it is as well to be aware of the ``hereditary'' version just stated.  We will in fact only make use of the following corollary, specific to the case where $H$ is a multiplicative subgroup of $\F_p^{\times}$. 

\begin{corollary}\label{add-spec-2}
Let $\alpha$, $0 < \alpha < 1$ be a parameter, let $H \subseteq \F^{\times}_p$ be a multiplicative subgroup and set $A := \Spec_{\alpha}(H)$, $A' := \Spec_{\alpha^2/2}(H)$. Write $L := |A'|/|A|$. 
Then for each $h \in H$ the additive energy $\omega_+(A,h\cdot A)$ is at least $\alpha^4/L$. 
\end{corollary}
\proof  For each $x \in \F_p$ write $r(x)$ for the number of pairs $a,a' \in A$ with $a - a' = x$. By Lemma \ref{add-spec} we know that $\sum_{x \in A'} r(x) \geq \frac{1}{2}\alpha^2|A|^2$. It follows from the Cauchy-Schwarz inequality that 
\[ \sum_{x} r(x)^2 \geq \sum_{x \in A'} r(x)^2 \geq \frac{1}{|A'|} \big( \sum_x r(x)\big)^2 = \frac{\alpha^4}{L}|A|^3.\] 
The left-hand side is the number of solutions to $a_1 + a_2 = a_3 + a_4$. However $A$ is $H$-invariant, and so this is equal to the number of solutions to $a_1 + ha_2 = a_3 + ha_4$.\endproof

The bound obtained here is of course strongest when $A' = \Spec_{\alpha^2/2}(H)$ is not a great deal larger than $A = \Spec_{\alpha}(H)$. The heart of the proof of Theorem \ref{bgk-theorem}, which we now present, consists of a pigeonholing argument which enables us to locate a scale $\alpha$ for which this happens.

Before doing that let us record the result of combining the last corollary with the additive-multiplicative Balog-Szemer\'edi-Gowers theorem.

\begin{proposition}\label{add-spec-4}
Suppose that $H$ is a multiplicative subgroup of $\F^{\times}_p$, $\alpha$ is a parameter with $0 < \alpha < 1$, $A = \Spec_{\alpha}(H)$, $A' = \Spec_{\alpha^2/2}(H)$ and $L = |A'|/|A|$. Using approximate notation at scale $L/\alpha$,  there is a set $X \subseteq A$, $|X| \approx |A|$, and a set $H' \subseteq H$, $|H'| \approx |H|$, such that  $|X + h \cdot X| \approx |X|$ for all $h \in H'$.\endproof
\end{proposition}

Suppose then that $H \leq \F_p^{\times}$ is a subgroup of size at least $p^{\delta}$ and that 
\[ |\sum_{x \in H} e(\xi x/p)| \geq \eta |H|\] for some $\xi\neq 0$, 
or in other words that $\Spec_{\eta}(H) \neq \{0\}$. Our aim is to obtain an upper bound of the form $\eta \ll p^{-\delta'}$ for some positive $\delta'$ (depending on $\delta$). Suppose for a contradiction that there is no such bound, thus $\eta = p^{-o(1)}$.

Let $J$ be an integer to be specified later, set $\alpha_0 := \eta$ and define $\alpha_{i+1} := \alpha_i^2/2$ for $i = 0,1,2,\dots,J$. We have the nesting
\[ \Spec_{\alpha_0}(H) \subseteq \Spec_{\alpha_1}(H) \subseteq \dots \subseteq \Spec_{\alpha_J}(H),\] and so by a trivial instance of the pigeonhole principle there is some $i$ such that 
\[ |\Spec_{\alpha_{i+1}}(H)| \leq p^{1/J}|\Spec_{\alpha_i}(H)|.\]
Set $\alpha := \alpha_i$, $A := \Spec_{\alpha_i}(H)$ and $A' := \Spec_{\alpha_{i+1}}(H)$, and note that $\alpha \geq (\eta/2)^{2^J} = p^{-o_J(1)}$. Moreover it follows from Parseval's identity that $|A| \leq p^{1 - \delta}\alpha^{-2}$, which is less than $p^{1 - \delta/2}$ for large $p$. Since $A$ is nonempty and $H$-invariant, we also have $|A| \geq p^{\delta}$.

Applying Proposition \ref{add-spec-4}, we obtain a set $X \subseteq A$ with $|X| \geq (\alpha/p^{1/J})^C|A|$ and a set $H' \subseteq H$ with $|H'| \geq (\alpha/p^{1/J})^C|H|$ such that $|X + h \cdot X| \leq (p^{1/J}/\alpha)^C|X|$ for all $h \in H'$. Using the bounds noted above, these inequalities may be written $|X| \geq p^{o_J(1) -C/J }|A|$, $|H'| \geq p^{o_J(1) - C/J}|H|$ and $|X + h \cdot X| \leq p^{C/J - o_J(1)}|X|$. In particular if $J$ is chosen large enough in terms of $\delta$ then we have $|X| \geq p^{\delta/2}$, $|H'| \geq p^{\delta/2}$ and $|X + h \cdot X| \leq p^{f(\delta)}|X|$, where $f(\delta)$ is any quantity depending on $\delta$ that might later prove convenient for us. By choosing $f(\delta)$ appropriately, we obtain a contradiction to Corollary \ref{corollary6.2}, thereby completing the proof of Theorem \ref{bgk-theorem}.\endproof

\section{Acknowledgements} It is a pleasure to thank Jean Bourgain and Elon Lindenstrauss for helpful communications.

\providecommand{\bysame}{\leavevmode\hbox to3em{\hrulefill}\thinspace}
\providecommand{\MR}{\relax\ifhmode\unskip\space\fi MR }
\providecommand{\MRhref}[2]{%
  \href{http://www.ams.org/mathscinet-getitem?mr=#1}{#2}
}
\providecommand{\href}[2]{#2}

     \end{document}